\documentclass[12pt]{article}
\usepackage{amsfonts}
\usepackage{eepic}
\usepackage{epic}
\usepackage{color}
\usepackage{graphicx}
\begin{document}
\newtheorem{Theoreme}{Th\'eor\`eme}[section]
\newtheorem{Theorem}{Theorem}[section]
\newtheorem{Th}{Th\'eor\`eme}[section]
\newtheorem{De}[Th]{D\'efinition}
\newtheorem{Pro}[Th]{Proposition}
\newtheorem{Lemma}[Theorem]{Lemma}
\newtheorem{Proposition}[Theorem]{Proposition}
\newtheorem{Lemme}[Theoreme]{Lemme}
\newtheorem{Corollaire}[Theoreme]{Corollaire}
\newtheorem{Consequence}[Theoreme]{Cons\'equence}
\newtheorem{Remarque1}[Theoreme]{Remarque}
\newtheorem{Convention}[Theoreme]{{\sc Convention}}
\newtheorem{PP}[Theoreme]{Propri\'et\'es}
\newtheorem{Conclusion}[Theoreme]{Conclusion}
\newtheorem{Ex}[Theoreme]{Exemple}
\newtheorem{Definition}[Theorem]{Definition}
\newtheorem{Remark1}[Theorem]{Remark}
\newtheorem{Not}[Theoreme]{Notation}
\newtheorem{Nota}[Theorem]{Notation}
\newtheorem{Propo}[Theorem]{Proposition}
\newtheorem{exercice1}[Th]{Lemme-Confi\'e au lecteur}
\newtheorem{Corollary}[Theoreme]{Corollary}
\newtheorem{PPtes}[Th]{Propri\'et\'es}
\newtheorem{Defi}[Theorem]{Definition}
\newtheorem{Example1}[Theorem]{Example}
\newenvironment{Proof}{\medbreak{\noindent\bf Proof }}{~{\hskip
3pt$\bullet$\bigbreak}}

\newenvironment{Demonstration}{\medbreak{\noindent\bf D\'emonstration
 }}{~{\hskip 3pt$\bullet$\bigbreak}} 

\newenvironment{Remarque}{\begin{Remarque1}\em}{\end{Remarque1}} 
\newenvironment{Remark}{\begin{Remark1}\em}{\end{Remark1}}
\newenvironment{Exemple}{\begin{Ex}\em}{~{\hskip
3pt$\bullet$}\end{Ex}} 
\newenvironment{exercice}{\begin{exercice1}\em}{\end{exercice1}}
\newenvironment{Notation}{\begin{Not}\em}{\end{Not}}
\newenvironment{Notation1}{\begin{Nota}\em}{\end{Nota}}
\newenvironment{Example}{\begin{Example1}\em}{~{\hskip
3pt$\bullet$}\end{Example1}}

\newenvironment{Remarques}{\begin{Remarque1}\em \ \\* }{\end{Remarque1}}
\newcommand{\Sp}{{\mathbb S}}
\renewcommand{\Re}{{\cal R}}
\renewcommand{\Im}{{\cal F}}
\newcommand{\finpreuve}{~{\hskip 3pt$\bullet$\bigbreak}}
\newcommand{\hp}{\hskip 3pt}
\newcommand{\hph}{\hskip 8pt}
\newcommand{\hphh}{\hskip 15pt}
\newcommand{\vp}{\vskip 3pt}
\newcommand{\vpv}{\vskip 15pt}
\newcommand{\IP}{{\mathbb{IP}}}
\newcommand{\rd}{{\mathbb R}^2}
\newcommand{\R}{{\mathbb R}}
\newcommand{\Hyper}{{\mathbb H}}
\newcommand{\Int}{{\mathbb I}}
\newcommand{\Boule}{{\mathbb B}(0,1)}
\newcommand{\Cantor}{{\mathbb K}}
\newcommand{\K}{{\mathbb K}}
\newcommand{\B}{{\mathbb B}}
\newcommand{\Z}{{\mathbb Z}}
\newcommand{\Nat}{{\mathbb N}}
\newcommand{\N}{{\mathbb N}}
\newcommand{\p}{{\mathbb P}}
\newcommand{\Esp}{{\mathbb E}}
\newcommand{\Complex}{{\mathbb C}}
\newcommand{\Ha}{{\cal H}}
\newcommand{\Harm}{{\bold H}}
\newcommand{\Lcal}{{\cal L}}
\newcommand{\ds}{\displaystyle}
\newcommand{\un}{\bold 1}
\newcommand{\Cone}{C(x,r,\epsilon ,\Phi)}
\newcommand{\Cn}{C(x,2^{-n},\epsilon ,\Phi )}
\newcommand{\Tranche}{W(x,r,\epsilon,\Phi)}
\newcommand{\Wn}{W(x,2^{-n},\epsilon,\Phi)}
\newcommand{\WFn}{W(x,2^{-n},\epsilon,\Phi)\cap F}
\newcommand{\ovec}{\overrightarrow}
\newcommand{\red}{{\bold R}}
\newcommand{\dimH}{\dim_{\Ha}}
\newcommand{\diam}{\mbox{diam}}
\newcommand{\diamit}{\mbox{\em diam}}
\newcommand{\para}{\vskip 2mm}
\newcommand{\cod}{\stackrel{\mbox{\tiny cod}}{\sim}}
\newcommand{\cardit}{\mbox{\em card}}
\newcommand{\card}{\mbox{card}}
\newcommand{\Sphere}{{\mathbb S}_d}
\newcommand{\dist}{\mbox{ dist}}
\newcommand{\distit}{\mbox{\em dist}}
\newcommand{\Tri}{{\cal P}}
\newcommand{\LL}{{\mathcal L}}
\newcommand{\infess}{\mbox{inf\,ess}}
\newcommand{\supess}{\mbox{sup\,ess}}

\definecolor{darkblue}{rgb}{0,0,.5}
\def\u{\underline}
\def\o{\overline}
\def\h{\hskip 3pt}
\def\hh{\hskip 8pt}
\def\hhh{\hskip 15pt}
\def\v{\vskip 8pt}
\def\vv{\vskip 15pt}
\font\courrier=cmr12
\font\grand=cmbxti10
\font\large=cmbx12
\font\largeplus=cmr17
\font\small=cmbx8
\font\nor=cmbxti10
\font\smaller=cmr8
\font\smallo=cmbxti10
\openup 0.3mm
\author{Athanasios BATAKIS and Viet Hung NGUYEN}
\title{On the exit distribution of partially reflected brownian motion in planar domains}
\maketitle

\begin{center}\begin{minipage}{15cm}
\noindent {\bf Abstract:} \em We show that the dimension of the exit distribution of planar partially reflected Brownian motion can be arbitrarily close to 2.
\end{minipage}
\end{center}
\section{Introduction} Let $\Omega$ be a domain in $\R^2$. It is well known (see \cite{Mak}, \cite{JW}) that the exit distribution of Brownian motion in $\Omega$ is carried by a borel subset of the boundary of dimension at most 1 (equal to one for simply connected domains). We are interested in the minimal dimension of sets carrying the exit distribution of partially reflected Brownian motion. 

The problem is posed as follows. Consider an $(\epsilon,\delta)$ domain $\Omega$ , take $F\subset \Omega$ a  closed subset of the boundary of $\Omega$ and consider Brownian Motion inside $\Omega$ absorbed by $F$ and reflected on $\partial \Omega\setminus F$ (for definitions of the $(\epsilon,\delta)$ domains and of reflected brownian motion see section \ref{definitions}). 
Note ${\mathcal R}_t$ the above process and $\tau_F$ the (first) hitting time of $F$ by ${\mathcal R}_t$. In general, $\tau_F$ may not be finite or may be finite but of infinite expectation (see also the so called ``trap domains'' \cite{BCM}).

We prove the following theorem


\begin{Theorem}\label{main}
For all $\eta>0$ there exist a domain $\Omega$ (that can be taken simply connected) and $F\subset\partial\Omega$ such that  $\p_x(\tau_F<\infty)=1$ and  for all $x\in\Omega$ and  for all $A\subset F$ of dimension $\dim A<2-\eta$ we have,
$$\p_x({\mathcal R}_{\tau_F}\in A)=0$$
\end{Theorem}

In particular this answers a question of B. Sapoval concerning Brownian motion as we will point out at the end of the paper: Consider a domain $\Omega$ and let $A$ be a subset of the boundary of (standard)  harmonic measure equal to $1$. If we change $A$ into reflecting boundary,  is the dimension of the exit distribution for this new diffusion still less than 1?

\noindent{\bf Acknowledgement:} The author would like to thank A. Ancona, L. Veron and M. Zinsmeister for many discussions that helped to clarify the original arguments and simplify the early proofs.
\section{Definitions of the main objects}\label{definitions}
The following definition is due to P. Jones \cite{Jones}.
\begin{Definition}\label{ed}
We say that a (not necessarily simply connected) domain $\Omega$ is an $(\epsilon,\delta)$-domain or locally uniform if there exist constants $\epsilon$ and $\delta$ such that for all $x,y\in \Omega$ with $|x-y|<\delta$ there is a (rectifiable) curve $\gamma$ joining $x$ and $y$ satisfying
\begin{enumerate}
\item $\epsilon\ell(\gamma)\le|x-y|$
\item $\epsilon\min\{|x-z|,|y-z|\}\le {\distit(z,\partial\Omega)}$
\end{enumerate}
\end{Definition}

 The $(\epsilon,\delta)$-domains satisfy the so called $W^{1,2}$- extension property, cf \cite{Jones}: if we note $W^{1,2}(\Omega)=\{f\in L^2(\Omega)\;;\;\nabla f\in L^2(\Omega)\}$ with the usual Sobolev norm $||f||_{1,2}=||f||_2+||\nabla f||_2$, we assume that there is a bounded linear operator $T: W^{1,2}(\Omega)\to W^{1,2}(\R)$ extending the identity of $ W^{1,2}(\Omega)$.

For $f,g\in W^{1,2}(\Omega)$ define 
$${\mathcal E}(f,g)=\int_{\Omega}<\nabla f,\nabla g>dx,$$ and 
$$ {\mathcal E}_1(f,g)={\mathcal E}(f,g)+\int_{\Omega} f gdx.$$

The Dirichlet form $({\mathcal E}, W^{1,2}(\Omega))$ is said to be regular on $\overline\Omega$ if
$W^{1,2}(\Omega)\cap C(\overline\Omega)$ is dense both in  $(W^{1,2}(\Omega),{\mathcal E}_1^{\frac12})$ and in 
$(C(\overline\Omega),||.||_{\infty})$.
Clearly, if $\Omega$ is a $(\epsilon,\delta)$-domain the Dirichlet form $({\mathcal E}, W^{1,2}(\Omega))$ is regular on $\overline \Omega$.

Following \cite{Chen}, \cite{BCR} we can now define the ``reflected'' Brownian motion. If $\Omega$ in an $(\epsilon,\delta)$-domain, there is a strong Markov process ${\mathcal R}$ associated, with continuous sample paths. Furthermore, we can construct a family of distributions $({\mathcal R}_t^x)_t$ for this process starting at every $x\in\overline\Omega$ (for further detail see also \cite{Fukushima}).

Take $F$ a closed subset of $\partial\Omega$ and consider $\tau_F$ the hitting time of $F$ for the process ${\mathcal R}_t^x$.
Now if we suppose that $\Esp_x[\tau_F]<+\infty$ for at least one $x\in \Omega$, we get that for any $f\in C(F)$, the function 
$$u:x\mapsto \Esp_x\left[f({\mathcal R}_{\tau_F})\right]$$ is bounded harmonic in $\Omega$ and takes the value $f$ at all regular points of $F$.

If we suppose that $\partial\Omega\setminus F$ is smooth then $u$ is the solution to the mixed Dirichlet-Neumann problem 
\begin{equation}\label{DN}
\left\{\begin{array}{l}
{\displaystyle u \mbox{ harmonic in }\Omega}\\
{\displaystyle\frac{\partial u}{\partial \eta}=0 \mbox{ on }\partial\Omega\setminus F}\\
{\displaystyle u=f\mbox{ on } F}
\end{array}\right.,
\end{equation}

where $\eta$ denotes the normal vector to the boundary $\partial\Omega$.

\begin{Remark}
We denote by $C_F(\overline\Omega)$ the set of continuous functions on $\overline\Omega$ vanishing on $F$. Suppose that
$ W^{1,2}(\Omega)\cap C_F(\overline\Omega)$ is dense in  $C_F(\overline\Omega),||.||_{\infty})$. We can then define the stochastic  process ${\mathcal R}^F_t$ associated. This process agrees with the previous one  for all $(\epsilon,\delta)$-domains  (see also \cite{AB}). 
\end{Remark}

Let $\omega_.$ denote the harmonic measure of this diffusion, ie. for $x\in\Omega$ and $A\subset \partial\Omega$, 
$$\omega_x(A)=\p_x({\mathcal R}_{\tau_F}\in A).$$

Remark  that, from (\ref{DN}), for $A\subset \partial \Omega$ measurable, the function $x\mapsto \omega_x(A)$ is positive harmonic in $\Omega$, tending to $1$ on $A$, to $0$ on $F\setminus A$ and of nul normal derivative on $\partial\Omega\setminus F$ . 

In the following we keep this same notation.
\section{Preliminary lemmas and remarks}
Let $E\subset\R^2$ be any set and, for every covering ${\mathcal V}_{\delta}(E)$ of $E$  with discs of radius less than $\delta$, let ${ H}_{\alpha}\left({\mathcal V}_{\delta}(E)\right)={\sum_{B\in {\mathcal V}_{\delta}} \diam(B)^{\alpha}}$. Consider 
$${\mathcal H}^{\delta}_{\alpha}(E)=\inf_{{\mathcal V}_{\delta}(E)}{ H}_{\alpha}\left({\mathcal V}_{\delta}(E)\right) \mbox{ and }
{\mathcal H}_{\alpha}(E)=\lim_{\delta\to 0}{\mathcal H}^{\delta}_{\alpha}(E)$$
Then, there exists an $\alpha_0\ge0$ such that ${\mathcal H}_{\alpha}(E)=0$  for all $\alpha>\alpha_0$ and ${\mathcal H}_{\alpha}(E)=\infty$ for all  $0\le \alpha<\alpha_0$. This $\alpha_0$ is denoted  $\dim_{\mathcal H}(E)$, the Hausdorff dimension of $E$.

For a Borel measure $\mu$ we define the Hausdorff dimension of $\mu$ as
$$\dim_{\mathcal H}(\mu)=\inf\{\dim_{\mathcal H}(E)\; ;\; \mu(E)>0\}$$

In particular, let $\mu$ be the harmonic measure $\omega_.$ defined above. Using the fact that, for any $A\subset F$,  $x\mapsto\omega_x(A)$ is harmonic we get that $\dim_{\mathcal H}(\omega_x)$ does not depend on the choice of $x\in\Omega$ and will be therefore denoted by $\dim_{\mathcal H}(\omega)$.

In this paper we are interested in the dimension of harmonic measure for partially reflected Brownian motion in domains in $\R^2$. Theorem \ref{main} can now be reformulated in the following terms:

{\em  ``for all $\eta>0$, there exists a uniform planar domain $\Omega$ and a closed set $F\subset \partial\Omega$  such that if $\omega_.$ is the harmonic measure for partially reflected Brownian motion (ie. reflected on $\partial\Omega\setminus F$, absorbed on $F$) we have $\dim_{\mathcal H}(\omega)>2-\eta$.''}

Clearly, $\dim_{\mathcal H}(\omega)\le\dim_{\mathcal H}(F)$. Therefore the boundary of $\Omega$ will be of Hausdorff dimension $\ge 2-\eta$.

\subsection{Potential theoretic lemmas}
By ``adapted cylinder"  ${\mathcal C}$ to a graph $\Gamma$ of a Lipschitz function $f$ we understand the intersection of a vertical revolution cylinder of finite height centered on $\Gamma$  with the $\Gamma^+=\{(x,y)\; ; \; y>f(x)\}$. We also ask  the ratio (height):(revolution radius) of  ${\mathcal C}$ to be greater than 2 times $||f||_L,$ la lipschitz norm of $f$.

We recall the boundary Harnack principle for reflected Brownian motion (see \cite{BassHsu}, \cite{Ancona1}). We say that $D$ is a Lipschitz domain if it is a Jordan domain and if the boundary is locally the graph of a Lipschitz function (with uniform lipschitz
norm).

\begin{center}
\begin{figure}
\begin{center}
\includegraphics[scale=0.3, width=9cm]{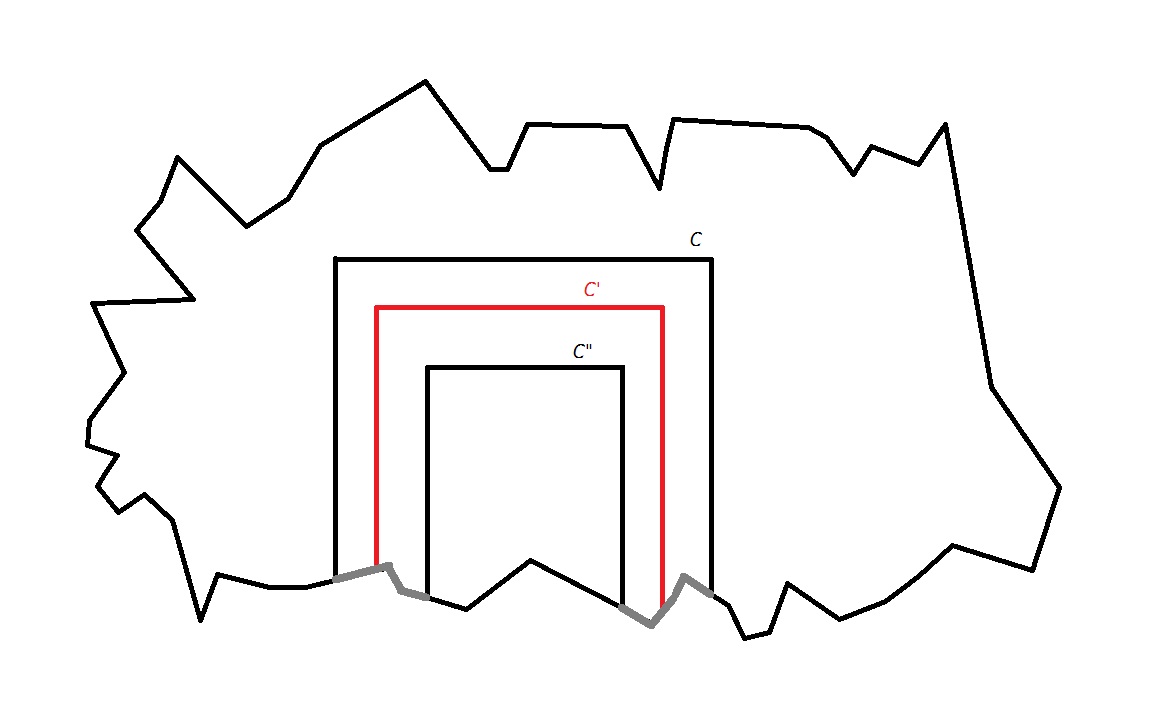}
\caption{\it Boundary Harnack Principle. \label{Harnack}}
\end{center}
\end{figure}
\end{center}

Let $D$ be a Lipschitz domain, $u$ and $v$ be positive harmonic functions on  $D$ with vanishing normal derivatives on the graph between the adapted cylinder  (to a graph-component of the boundary) ${\mathcal C}$ and  the ``sub"-adapted cylinder ${\mathcal C}''$ of the same center and revolution axis but of $\ell$ times the size, $\ell<1$ (see figure \ref{Harnack}).
\begin{Proposition}\label{harnack1}
 If ${\mathcal C}'$  is the ``middle" cylinder of the same center and revolution axis but of $\frac{1+\ell}{2}$ times the size of ${\mathcal C}$. Then for all   $x\in\partial  {\mathcal C}'\cap V$
$$\frac{v(x)}{u(x)}\sim \frac{v(P)}{u(P)},$$
where $P$ is the intersection point of the revolution axis of the cylinder ${\mathcal C}'$ and of its boundary.
\end{Proposition}
 The multiplicative constants in the equivalence relation depend on the ratio (revolution radius):(height) of ${\mathcal C}$ , on $\ell$, on the Lipschitz norm of the boundary and on the dimension of the space $n$ (here $n=2$) see also \cite{Ancona2}.

We also need a Dirichlet-Neumann version of the maximum principle. 
\begin{Proposition}\label{max1}
Let $D$ be a planar  domain , $\Gamma$ a continuous subset of the boundary of $D$, graph of a Lipschitz function, and $u$ a function harmonic in $D$ such that 
$\liminf_{y\to x} u(y)\ge 0$ for all $x\in\partial D\setminus \Gamma$ and $\frac{\partial u}{\partial\eta }=0$ on $\Gamma$, where $\eta$ denotes the  normal vector on $\Gamma$.
Then $u\ge 0$ on $D$.
\end{Proposition}
This is a consequence of the unicity of solutions (see for instance \cite{Horm}) and the probabilistic description of these same solutions of the mixed Dirichlet-Neumann problem, described above.

\subsection{Subsidiary results}
We will use the following result due to Benjamini, Chen and Rohde.
\begin{Theorem}(Theorem 5.1 of \cite{BCR}) \label{BCR}
 Let $\Omega$ be a locally uniform bounded planar domain. Then, $\dim_{\mathcal H}\left( R\left([0,\infty)\right)\cap\partial\Omega\right)=\dim_{\mathcal H}(\partial\Omega)$, $\p_x$-almost surely, for all $x\in\overline\Omega$.
\end{Theorem}
In particular, under the assumptions of the theorem, if $F\subset\Omega$ is a closed set such that $\dim_{\mathcal H}(\partial\Omega\setminus F)<\dim_{\mathcal H}(F)$ we have 
\begin{equation}\label{finitepro}
\p_x\left (\tau_F<+\infty\right)=1
\end{equation}
for all $x\in\overline\Omega$.

\begin{Proposition}\label{4.2}
Under the same assumptions, formula (\ref{finitepro}) implies 
$\Esp_x\left[\tau_F\right]<+\infty,$
for all $x\in \overline \Omega$.
\end{Proposition}
To prove this proposition we recall a result of Burdzy, Chen and Marshall.
\begin{Theorem}(\cite{BCM}) \label{BCM}
Let $\Omega$ be any bounded locally uniform domain and $\B$ a closed ball in $\Omega$. If we note $\tau_{\B}$ the hitting time of 
$\B$ by ${\mathcal R}_t$ then 
$\sup_{x\in\overline\Omega}\Esp_x\left[\tau_{\B}\right]<\infty$.
\end{Theorem}

\begin{Proof}{\bf of Proposition \ref{4.2}. }
Let $z$ be a point in $\Omega$ and $\B_z\subset\Omega$ a closed disc centered at $z$. For $x\in\Omega$ and any $s>0$,
$$\Esp_x\left[\tau_F\right]\le \sum_{n\in\N}s\p_x(\tau_F\ge ns).$$
By formula (\ref{finitepro}) for all $N\in\N$ there exists $s$ sufficiently big such that 
$\p_x\left (\tau_F>\frac{s}{2}\right)<\frac1N$ and $\p_z\left (\tau_F>\frac{s}{2}\right)<\frac1N$. Furthermore we can choose $s>2\sup_{x\in\overline\Omega}\Esp_x\left[\tau_{\B_z}\right]$.
We get that, for $n>1$,
$$\p_x(\tau_F\ge ns)=\p_x(\tau_F\ge ns|\tau_F\ge (n-1)s) \p_x(\tau_F\ge (n-1)s)$$
We can bound $\p_x(\tau_F\ge ns|\tau_F\ge (n-1)s)\le \sup_{y\in \Omega}\p_y\left(\tau_F>s\right)$.
On the other hand $\p_y\left(\tau_F>s\right)\le \p_y\left(\tau_{\B_z}>s/2\right)+\p_y\left(\tau_{\B_z}<s/2\; , \;{\mathcal R}_{[\tau_{\B_z},s]}\cap F=\emptyset\right).$

Using the Markov property of ${\mathcal R}$ , 
$$\p_y\left(\tau_{\B_z}<s/2\; , \;{\mathcal R}_{[\tau_{\B_z},s]}\cap F=\emptyset\right)\le \p_y\left(\tau_{\B_z}<s/2\right) \sup_{v\in\B_z}\p_v\left(\tau_F>s/2\right).$$
Using parabolic Harnack principle (see \cite{BCM}) we get that there is a constant $c>1$ such that 
$$\sup_{v\in\B_z}\p_v\left(\tau_F>\frac{s}{2}\right)\le c \p_z\left (\tau_F>\frac{s}{2}\right)<c/N.$$
We also have 
$$ \p_y\left(\tau_{\B_z}<s/2\right)<\frac{2\sup_{x\in\overline\Omega}\Esp_x\left[\tau_{\B_z}\right]}{s}$$ therefore, for $s$ big enough,
$$\p_x(\tau_F\ge ns|\tau_F\ge (n-1)s)< 1/2$$
By induction we get $\p_x(\tau_F\ge ns)\le \left({\frac12}\right)^n$ and hence $\Esp_x\left[\tau_F\right]<+\infty$.
\end{Proof}
In fact we have proved that $\sup_{x\in\Omega}\Esp_x\left[\tau_F\right]<+\infty$.
\section{Proof of theorem \ref{main}}
Even though our proof can be carried out using only simply connected domains we have chosen to present a totally disconneted example: the constructions  appear better and the lemmas get easier to write.
\subsection{Construction of the domain}

\begin{center}
\begin{figure}
\begin{center}
\includegraphics[scale=0.3, width=9cm]{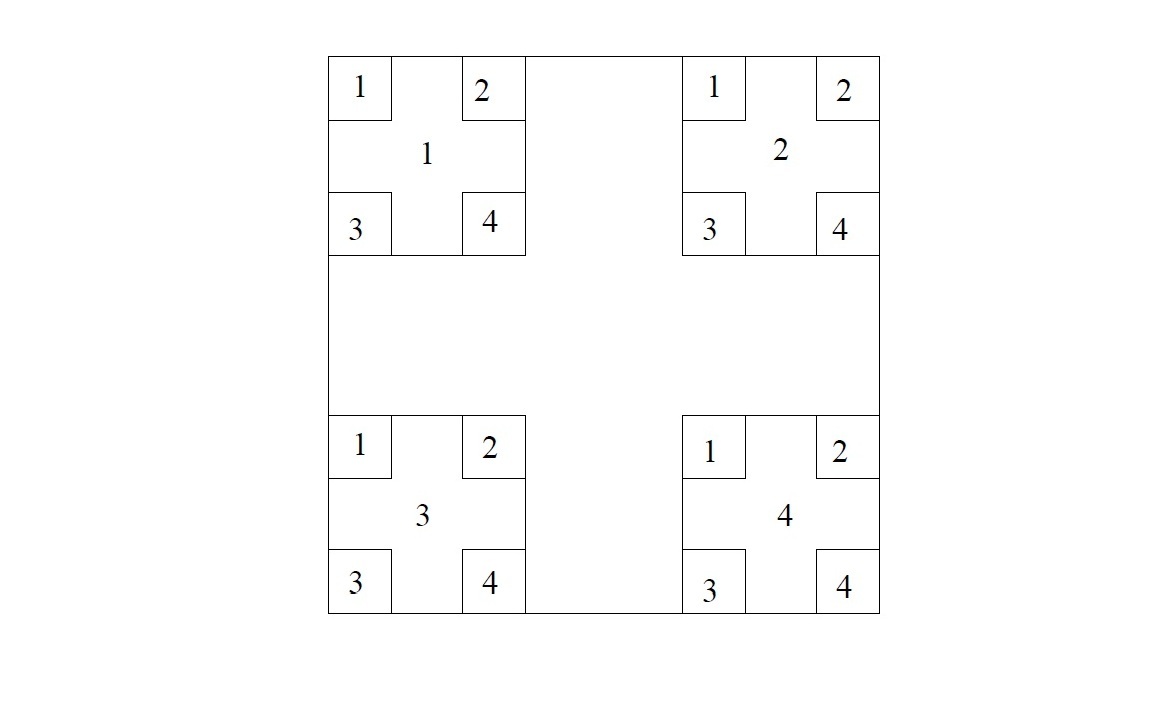}
\includegraphics[scale=1, width=6cm]{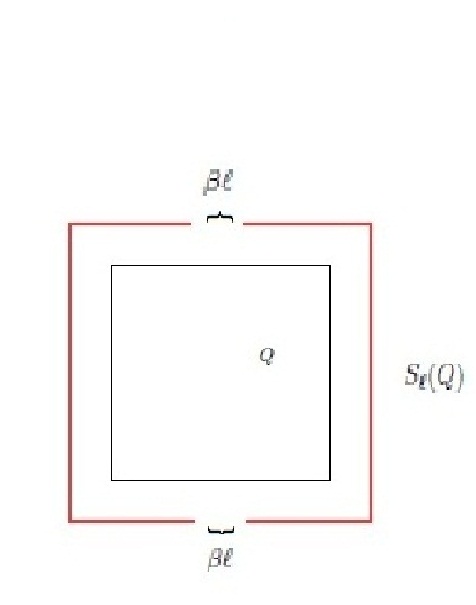}
\caption{\it  A. $4$-corner Cantor set and its encoding. \label{cantor}  B. The squares $S_{\ell} (Q)$ . \label{squares}}
\end{center}
\end{figure}
\end{center}

We construct, for $\alpha\in(0,\frac12)$ a $4$-corner Cantor set (fig \ref{cantor}.A) in the following way. We start with the square $Q=[-\frac12,\frac12]^2$ that we replace by four squares of sidelength $\alpha$ situated at the four corners of $Q$.  We name these squares $Q_1,...,Q_4$. We replace then each $Q_i$ , $i=1,...4$ by four smaller squares of sidelength $\alpha^2$ situated at the corners of $Q_i$. We note these squares of the second generation $Q_{ij}$, where $j=1,...,4$ and so on.
Let us denote $\K$ the Cantor set constructed in this way. We endowe $\K$ with the natural encoding  identifying it to the abstract Cantor set $\{1,...,4\}^{\N}$.

Observe that $\dim_{\mathcal H}(\K)=\left|\frac{\log4}{\log{\alpha}}\right|$ and hence for $\alpha$ close to $\frac12$ the dimension of the Cantor set is close to $2$.
The set $\K$ will be the absorbing part of the boundary of $\Omega$. Let us know construct the reflecting part. First of all, in order to ensure boundedness let us consider a ball $\B_0$, centered at $0$ of radius, say, $10^6$. The domain $\Omega$ will be a subset of $\B_0\setminus \K$.

Let $Q$ be a square of sidelength $\rho$ centered at $(x^*,y^*)$ and for $0<\beta<10^{-2}\rho$ and $\ell>1$ consider the ``unfinished'' squares 
$$S_{\beta,\ell} (Q)=\{(x,y)\in\R^2\;;\; |x-x^*|=|y-y^*|=\rho \ell  \; \mbox{and}\;x\notin (x^*-\beta\ell/2,x^*+\beta\ell/2) \}$$ (see figure \ref{cantor}.B).
Finally consider the blown-up version of $S_{\beta,\ell} (Q)$ (see figure \ref{Foot}):
\begin{equation}\label{Ls}
L_{\beta,\ell} (Q)=\{z\in\R^2\;;\;\dist(z, S_{\beta,\ell} (Q))\le 10^{-6}\beta\rho\}
\end{equation}

Note that, if $\ell$ is less than $\frac{1}{2\alpha}$, for any $Q$ and $Q'$  squares of the construction of the Cantor set  $L_{\beta,\ell} (Q)\cap L_{\beta,\ell} (Q')=\emptyset$. 
Consider the union of the Cantor set $\K$ with $$M_{\beta,\ell}(\K)=\displaystyle \bigcup_{n\in\N}\bigcup_{i_1,...i_n}\left(L_{\beta,\ell} \left(Q_{i_1,...i_n}\right)\right)$$
The domain $\Omega$ is defined as the complementary of this union within the ball $\B_0$ of radius $10^6$:
$$\Omega=\B_0\setminus\left(\K\cup M_{\beta,\ell}(\K)\right).$$

\begin{Remark}
For $\ell>1$ and $\beta>0$ fixed the domain $\Omega$ is clearly a bounded uniform domain. Therefore we can construct partially reflected Brownian motion ${\mathcal R}_t$ in $\Omega$ with the partition of the boundary of $\Omega$ into an absorbing part of the boundary $F=\K$ and a reflecting part  $\partial\Omega\setminus F$. 

Note also that the Hausdorff dimension of $\partial\Omega $ equals the Hausdorff dimension of $\K$ if $\dim_{\mathcal H}\K>1$ (ie. if $\alpha>\frac 14$). This is because $\partial\Omega\setminus K$ consists of a countable union of rectifiable arcs and is therefore of Hausdorff dimension $1$.
\end{Remark}

We will show that for $\alpha\in (0,\frac12)$ and every $\epsilon>0$ there exists  $\ell<\frac{1}{2\alpha}$ and $\beta$ close to $0$ such that the domain $\Omega$, constructed in the previous way, satisfy $\dim_{\mathcal H}\omega>(1-\epsilon)\dim_{\mathcal H}\K$.

\subsection{Preparatory lemmas}

\begin{figure}
\begin{center}
\includegraphics[scale=1, width=10cm]{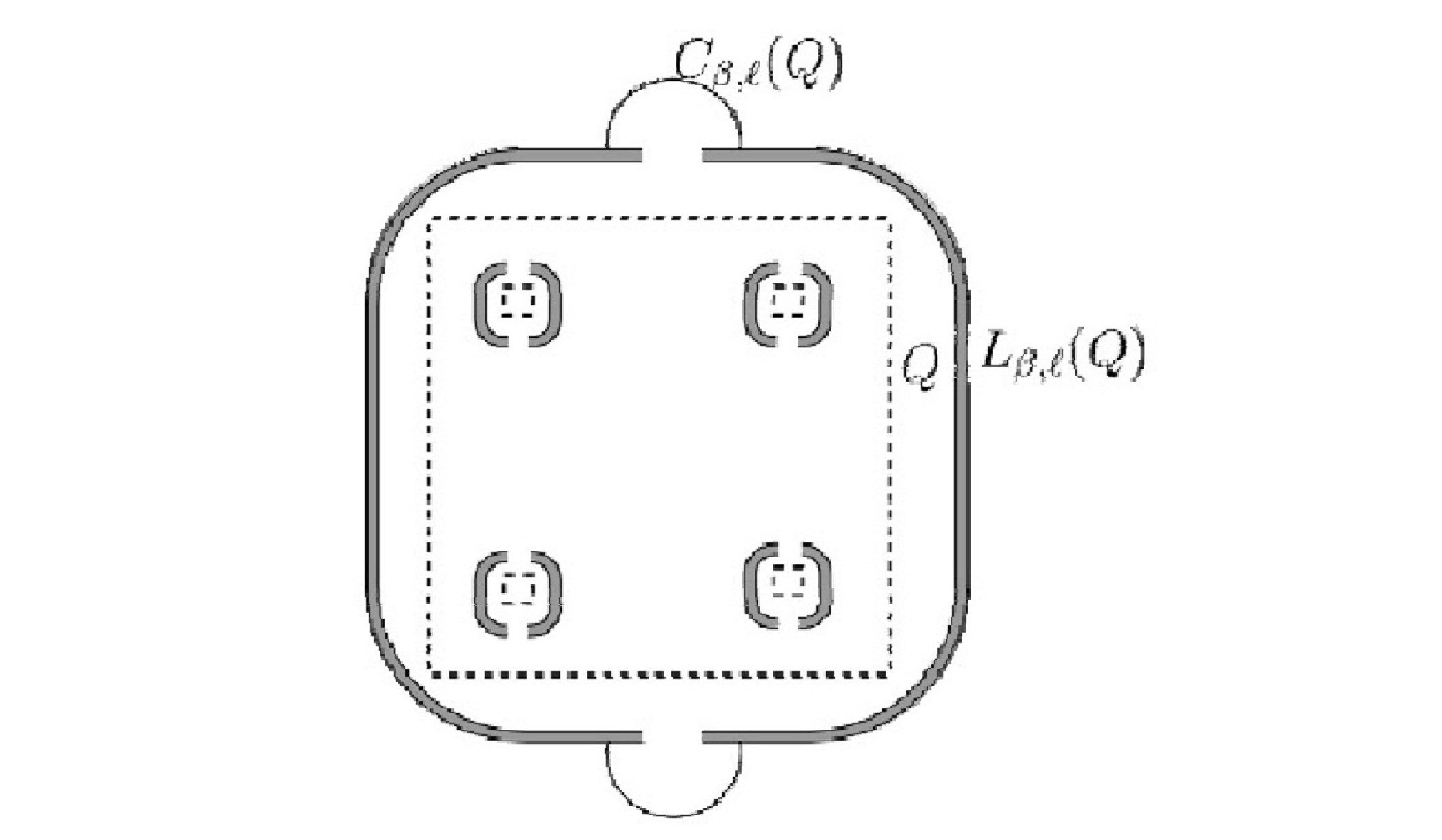}
\caption{$F_{\beta,\ell}(Q)$\label{Foot} }
\end{center}
\end{figure}

Let  $Q=Q_{i_1,...,i_n}$ be a square of the construction of $\K$ of sidelength $\rho$  and let $(x^*,y^*)$ be it's center. Let $F_{\beta,\ell}(Q)=L_{\beta,\ell} (Q)\cup C_{\beta,\ell}(Q)$ where 
\begin{eqnarray*}
C_{\beta,\ell}(Q)&=&\{(x,y)\;;\; y>y^*+\rho\ell\mbox{ and } ||(x-x^*,y-y^*-\rho\ell)||=(\frac{1}{2\alpha}-\ell)\rho\}\\
&&\bigcup \{(x,y)\;;\; y<y^*-\rho\ell\mbox{ and } ||(x-x^*,y-y^*+\rho\ell)||=(\frac{1}{2\alpha}-\ell)\rho\},
\end{eqnarray*}
see figure \ref{Foot}.

Let $\tilde D$ be the bounded component of the complementary of  $F_{\beta,\ell}(Q)$ and $D=\Omega\cap \tilde D$ (as in figure \ref{Foot}). 

We consider reflected Brownian motion ${}^D{\mathcal R}$ in $D$ and we note $\tau^D$ the hitting time of $\K\cup C_{\beta,\ell}(Q)$ by ${}^D{\mathcal R}$.
It follows on the previous discussion that for all $x\in D$,  $\Esp_x{\tau^D}<\infty$ and, furthermore, $\p_x\left({}^D{\mathcal R}_{\tau^D}\in \K\right)>0$. To prove this last claim one can also use the arguments of relation (\ref{V_i}) below, applied to the domaine $D$ and to the diffusion ${}^D{\mathcal R}$ respectively.

\begin{Remark}
Consider $Q_1=Q_{i_1,...,i_na},Q_2=Q_{i_1,...,i_nb}$ (with $a,b=1,..,4$) two sub-cubes of  $Q=Q_{i_1,...,i_n}$. By symmetry we get that, if $x_Q=(x^*,y^*)$ is the center of $Q$,
$$\p_{x_Q}\left({}^D{\mathcal R}_{\tau^D}\in \K\cap Q_1\right)=\p_{x_Q}\left({}^D{\mathcal R}_{\tau^D}\in \K\cap Q_2\right).$$
It follows, using Harnack's principle, that for all $\epsilon>0$ there exists $r=r_{\epsilon}>0$ (depending only on $\epsilon,\ell,\alpha$ but not on $\beta$)  such that 
\begin{equation}\label{epsilon}
||x-x_Q||<r\Longrightarrow  \p_{x_Q}\left({}^D{\mathcal R}_{\tau^D}\in \K\cap Q_1\right)\le(1+\epsilon)\p_{x_Q}\left({}^D{\mathcal R}_{\tau^D}\in \K\cap Q_2\right).
\end{equation}
\end{Remark}

Let us now prove that for $\beta$ small enough, the harmonic  functions (measures) $U_i(.)=\p_{.}\left({}^D{\mathcal R}_{\tau^D}\in \K\cap Q_i\right)$, $i=1,2$,  satisfy inequality  (\ref{epsilon}) in the subdomain $D'$ of $D$ :

\begin{equation}\label{far}
D'=D\setminus\left(\B\left((x^*,y^*+\rho\ell),(\frac{1}{2\alpha}-\ell)\rho\right)\cup\B\left((x^*,y^*-\rho\ell),(\frac{1}{2\alpha}-\ell)\rho\right))\cup\bigcup_{i=1...4}\ell Q_i\right).
\end{equation} 

Remark that the closure of $D'$ is a compact subset of $D\cup L_{\beta,\ell} (Q)$.
\begin{Lemma}\label{zeta}
For all $\epsilon>0$ there exists $\beta_0>0$ such that for all $\beta<\beta_0$ and all  $x\in D'$
$$\p_x\left(\exists t<\tau^D\; ; {}^D{\mathcal R}_t\in\B(x_Q,r_{\epsilon})\right)>1-\epsilon$$
\end{Lemma}

\begin{Proof}
We introduce an auxilliary subdomain $D''$ of $D$, $D''=D\setminus \overline B(x_Q,r_{\epsilon})$. Consider, in $D''$ the harmonic function $\zeta$ satisfying the mixed Dirichlet-Neumann boundary conditions
$ \zeta=0$ on $B(x_Q,r_{\epsilon})$, 
$\zeta=1$ on $\K\cup C_{b,\ell}(Q)$
and 
$\frac{\partial \zeta}{\partial \eta}=0$ elsewhere on $\partial D''$.

It is immediate that, since $r_{\epsilon}$ does not depend on $\beta$, $\zeta$ tends to $0$  when $\beta$ goes to $0$. By the maximum principle \ref{max1}, 
$$1-\zeta(x)<\p_x\left(\exists t<\tau^D\; ; {}^D{\mathcal R}_t\in\B(x_Q,r_{\epsilon})\right).$$
On the other hand, for every $x$ there is an $\beta_0$ such that $\zeta(x)<\epsilon$ for all $\beta<\beta_0$ and by Harnack's principle this inequality can be taken uniform in $\overline{ D'}$.
\end{Proof}
Keeping the same notation we also have: 
\begin{Lemma}\label{zeta1}
Let $Q_i=Q_{i_1,...,i_ni}$, $i=1,...4$,  be a sub-cube of  $Q=Q_{i_1,...,i_n}$. Then, there is a constant $C>0$ depending only on $\alpha,\ell$ such that  for all $x\in C_{\beta,\ell}(Q_i)$ 
$$\p_x(\exists\; \;  0<t_1<t_2<\tau^D\; ;{}^D{\mathcal R}_{t_1}\in B(x_Q,r_{\epsilon})\; ,\;{}^D{\mathcal R}_{t_2}\in  C_{\beta,\ell}(Q_i))\ge C$$
\end{Lemma}
The proof of this lemma is standard and hence omitted.

\subsection{Harmonic measure estimates}
As before, let $Q_1=Q_{i_1,...,i_na},Q_2=Q_{i_1,...,i_nb}$ be two sub-cubes of a given cube $Q=Q_{i_1,...,i_n}$ of the construction of $\K$.

Take $U_1$ and $U_2$ to be the harmonic functions previously defined in $D$, ie. satisfying the mixed Dirichlet-Neumann boundary conditions :
\begin{equation}
\left\{\begin{array}{l} 
 U_1=1\mbox{ on }\K\cap Q_1\\
U_1=0\mbox{ on }(\K\cap Q_2)\cup C_{\beta,\ell}(Q)\\
\frac{\partial U_1}{\partial\eta}=0 \mbox{ elsewhere on }\partial D \\
\end{array}\right. \mbox{ and }
\left\{\begin{array}{l} 
U_2=1\mbox{ on }\K\cap Q_2\\
U_2=0\mbox{ on }(\K\cap Q_1)\cup C_{\beta,\ell}(Q)\\
\frac{\partial U_2}{\partial\eta}=0 \mbox{  elsewhere on  }\partial D\\
\end{array}\right.
\end{equation}
Thus, $U_i(.)=\p_{.}\left({}^D{\mathcal R}_{\tau^D}\in \K\cap Q_i\right)$, $i=1,2$.
\begin{Lemma}\label{hmD}
For every $\epsilon>0$ there existe $\beta_0>0$ such that for all  $0<\beta<\beta_0$ and all $x\in D'$,
$$ U_1(x)<(1+\epsilon) U_2(x).$$
\end{Lemma}

\begin{Proof}
The proof relies on lemma \ref{zeta}. By Harnack's principle there exists $C>0$ such that $U_i(x)\ge C U_i(x_Q)$,  for all $x\in D'$. On the other hand,
\begin{eqnarray*}
U_i(x)&=& \p_{.}\left({}^D{\mathcal R}_{\tau^D}\in \K\cap Q_i\right)=\\
&=& \p_{.}\left({}^D{\mathcal R}_{\tau^D}\in \K\cap Q_i\mbox{ , } {}^D{\mathcal R}_{[0,\tau^D]}\cap\B(x_Q,r_{\epsilon})=\emptyset\right)+\\
&+&\p_{.}\left({}^D{\mathcal R}_{\tau^D}\in \K\cap Q_i\mbox{ , } {}^D{\mathcal R}_{[0,\tau^D]}\cap\B(x_Q,r_{\epsilon}) \not=\emptyset\right)
\end{eqnarray*}
 
Lemmas \ref{zeta} and  \ref{zeta1} imply that $$ \p_{.}\left({}^D{\mathcal R}_{\tau^D}\in  Q_i\mbox{ ,} {}^D{\mathcal R}_{[0,\tau^D]}\cap\B(x_Q,r_{\epsilon})=\emptyset\right)\le \epsilon\p_{.}\left({}^D{\mathcal R}_{\tau^D}\in Q_i\mbox{ ,} {}^D{\mathcal R}_{[0,\tau^D]}\cap\B(x_Q,r_{\epsilon}) \not=\emptyset\right)$$
On the other hand,  by the Markov property,
$$\p_{.}\left({}^D{\mathcal R}_{\tau^D}\in Q_i\mbox{ ,} {}^D{\mathcal R}_{[0,\tau^D]}\cap\B(x_Q,r_{\epsilon}) \not=\emptyset\right)\le \sup_{x\in\B(x_Q,r_{\epsilon})}\p_{x}\left({}^D{\mathcal R}_{\tau^D}\in Q_i\right).$$
Therefore, using once more Harnack's inequality  
$$U_i(x)\le (1+c\epsilon)U_i(x_Q),$$ and the lemma's claim follows using symmetry.
\end{Proof}

Consider now the functions $V_1$ and $V_2$ that solve the following mixed Dirichlet-Neumann problem in $\Omega$.
\begin{equation}
\left\{\begin{array}{l} 
V_1=0\mbox{ on }\K\cap Q_1\\
V_1=U_1\mbox{ on }\K\setminus Q_1\\
\frac{\partial V_1}{\partial\eta}=\frac{\partial U_1}{\partial\eta} \mbox{ on }\partial\Omega\setminus \K\\
\end{array}\right. \mbox{ and }
\left\{\begin{array}{l} 
V_2=0\mbox{ on }\K\cap Q_2\\
V_2=V_1\mbox{ on }\K\setminus Q_2\\
\frac{\partial V_2}{\partial\eta}=\frac{\partial U_2}{\partial\eta} \mbox{ on }\partial\Omega\setminus \K\\
\end{array}\right.
\end{equation}
\subsection{Proof of theorem}
We need to show that for $\beta$ small enough $V_i\le (1+\epsilon)V_j$, for $i,j=1,2$. Since  $V_i(x)=\p_x({\mathcal R}_{\tau_{\K}}\in \K\cap Q_i)$, this inequality clearly implies that the harmonic measure for partially reflected Brownian motion $\omega$ satisfies $$(1+\epsilon)^{-n}4^{-n}\le\omega(Q_{i_1,...,i_n})\le (1+\epsilon)^n4^{-n},$$ for all $n$ and all indices $i_1,...i_n\in\{1,...,4\}$ and hence the claim.

 We note 
$$C'_{\beta,\ell}=D\cap\partial\left(\B\left((x^*,y^*+\rho\ell),(\frac{1}{2\alpha}-\ell)\rho\right)\cup\B\left((x^*,y^*-\rho\ell),(\frac{1}{2\alpha}-\ell)\rho\right)\right)$$
For $n\in\N$, consider the increasing sequences of hitting times
$$T_n=\inf\{t\; ;\; \exists t_1<s_1...<t_{n-1}<s_{n-1}<t\mbox{ s.t. } {\mathcal R}_{t_i}\in C'_{\beta,\ell}\; ,\; {\mathcal R}_{s_i}\in C_{\beta,\ell}\}$$ and $$S_n=\inf\{s\; ;\; \exists t_1<s_1...<t_n<s\mbox{ s.t. } {\mathcal R}_{t_i}\in C'_{\beta,\ell}\; ,\; {\mathcal R}_{s_i}\in C_{\beta,\ell}\}$$
with the convention $T_n=\infty$ ($S_n=\infty$) if the corresponding set is empty.
\begin{eqnarray}\label{V_i}
V_i(x)&=&\p_x({\mathcal R}_{\tau_{\K}}\in Q\cap\K\; , \;{\mathcal R}_{\tau_{\K}}\in \K\cap Q_i)\nonumber\\
&=&\sum_n\p_x\left(0<T_1<...<T_n<\tau_{\K}<\infty\;, \,S_n=\infty\; , \;{\mathcal R}_{\tau_{\K}}\in \K\cap Q_i\right)\nonumber\\
&=&\sum_n\p_x\left(0<T_1<...<T_n<\infty\right)\Esp_x\p_{R_{T_n}}\left(S_1=\infty\; , \; {\mathcal R}_{\tau_{\K}}\in \K\cap Q_i\right)
\end{eqnarray}
where the last equality is derived by Markov's property.
By lemma \ref{hmD}, for $i,j=1,2$ and all $z\in C'_{\beta,\ell}$,
$$\p_z\left(S_1=\infty\; , \; \tau_{\K}\in \K\cap Q_i\right)\le (1+\epsilon)\p_z\left(S_1=\infty\; , \; {\mathcal R}_{\tau_{\K}}\in \K\cap Q_j\right).$$
Implementing this inequality in (\ref{V_i}) we get 
$$V_i(x)\le \sum_n\p_x\left(0<T_1<...<T_n<\infty\right)\Esp_x\p_{R_{T_n}}\left(S_1=\infty\; , \;{\mathcal R}_{\tau_{\K}}\in \K\cap Q_j\right)$$ and hence
$V_i\le (1+\epsilon) V_j$, which completes the proof.

{\bf Comments-Further Remarks: }
If the boundary of  the domain $\Omega$ is entirely absorbing, ie. for the Laplace equation  with Dirichlet boundary conditions, then harmonic measure is carried by $\partial\Omega\setminus \K$. This is not difficult to see. In fact, using previous notation, to get to $\K$ Brownian motion has to go through an infinity of conformal annuli of the type $L{\beta,\ell}(Q)\setminus Q$. But, at every passage of this type, there is a -bounded from below probability- to hit $L{\beta,\ell}(Q)$. Hence the probability that brownian motion hits $\K$ is $0$.
Moreover, harmonic measure will be carried by a union of curves of finite length.

This answers the question of B. Sapoval mentionned in the introduction.

Another way to study the passage between  Dirichlet boundary condition to the Neumann boundary condition through the mixed Dirichlet-Neumann is through a random approach. This the object of a previous work \cite{BLZ} that should be completed in a forthcoming artcile \cite{ABZ}. 

\bibliographystyle{alpha}
\bibliography{biblio}

\end{document}